\documentclass{amsart}
\usepackage{amsmath}
\usepackage{amssymb}
\usepackage{latexsym}
\usepackage[dvips]{graphics}

\author{Ron Aharoni}
\address{Department of Mathematics\\Technion, Haifa\\ Israel 32000}
\email[Ron Aharoni]{ra@tx.technion.ac.il}
\thanks{The research of the first author was 
supported by the fund for the promotion of research at the Technion} 

\author{Eli Berger}
\email[Eli Berger]{beli@ibm.il.com}

\author{Ran Ziv}
\address{Department of Computer Science, Tel-Hai college, Upper Galilee,
Israel 12210}
\email[Ran Ziv]{ranziv@telhai.ac.il}

\title{A tree version of K\"onig's theorem}
\begin{document}
\maketitle

\newtheorem{theorem}{Theorem}[section]

\newtheorem{corollary}[theorem]{Corollary}
\newtheorem{notation}[theorem]{Notation}
\newtheorem{proposition}[theorem]{Proposition}
\newtheorem{assertion}[theorem]{Assertion}
\newtheorem{lemma}[theorem]{Lemma}
\newtheorem{definition}[theorem]{Definition}
\newtheorem{conjecture}[theorem]{Conjecture}

\newtheorem{remark}[theorem]{Remark}
\newtheorem{observation}[theorem]{Observation}

\newcommand{\gs}{\mbox{$G^{*k}$}}
\newcommand{\rtwo}{\mbox{$\Re^2$}}
\newcommand{\enp}{\mbox{$\hfill \Box$}}
\newcommand{\pf}{{\em Proof~}}
\newcommand{\calf}{\mbox{$\mathcal F$}}
\newcommand{\ger}{\mbox{$\succ$}}
\newcommand{\gep}{\mbox{$\sqsupset$}}

\begin{abstract}

K\"onig's theorem states that the covering number and the matching number of
a bipartite graph are equal. We prove a generalisation,
in which the point in one fixed side of the graph of each edge 
is replaced by a subtree of a given tree. 
The proof uses 
a recent extension of Hall's theorem to families of hypergraphs, by 
the first author and P. Haxell \cite{ah}.
\end{abstract}

\section{Width and matching width of families of trees}

A {\em matching} in a hypergraph is  a set of disjoint edges.
In \cite{ah} (see also \cite{ac} and \cite{meshulam}) there were introduced 
a few notions which proved to be fruitful in the study of 
matchings. The two main ones are those of {\em width} and {\em matching
  width}, which are presented below.

Let $H$ and $F$ be two hypergraphs on the same vertex set. 
A subset $C$ of  $F$ is said to be $H$-{\em covering} 
(or ``an $H$-{\em covering}'') 
if every edge in $H$ meets some edge from $C$ 
(in other words, the union of $C$ is a vertex cover for $H$). 
The {\em $F$-width} $w(H,F)$ of  $H$ is the minimal 
size of an $H$-covering set of edges from $F$.
%The {\em width} $w(H)$ of $H$ is $w(H,H)$.
The {\em $F$-matching width} $mw(H,F)$ 
of $H$ is the maximum, over all matchings $M$
in $H$, of $w(M,H)$.  We write  $w(H)$ and $mw(H)$
for $w(H,H)$ and $mw(H,H)$, respectively,
and call these parameters the {\em width} and {\em matching width} 
of $H$. 
These  notions have counterparts in which the covering set is 
required to be a matching.
The 
{\em independent  $F$-width} $iw(H,F)$ of $H$ is the minimal size of a  
covering matching in $F$. 
The {\em independent matching width} $imw(H,F)$ of $H$
is the maximum, over all matchings $M$ in $H$, of $iw(M,F)$.
Again, we write $iw(H)$ for $iw(H,H)$ and 
$imw(H)$ for $imw(H,H)$, and call them the {\em independent width} and
{\em independent matching width} of $H$, respectively. 
(The notion of independent matching width is used in \cite{ah},
but is not given there a special name.)

In \cite{ah} the following theorem was proved:

\begin{theorem}
\label{ah}
If $\mathcal A$ is a family of hypergraphs such that 
$mw(\bigcup {\mathcal B})\ge |{\mathcal B}|$
for every $\mathcal B \subseteq \mathcal A$, then there
is a choice of disjoint edges, one from each hypergraph in $\mathcal A$.
\end{theorem}

Clearly, $mw \le w$ and $imw \le iw$, and in general strict inequalities 
occur.
However,
as we shall prove below, for families of trees 
equality holds in both. In fact, both equalities
are special cases of a more general result.
For two hypergraphs $H_1, H_2$ and a symmetric relation $\sim$
on $H_2$ we denote by $w(H_1,H_2,~\sim)$
the minimal size of a $\sim$-related subset of $H_2$
which is an  $H_1$-covering.
By $mw(H_1,H_2,~\sim)$ we denote the maximum, over all matchings $M$ in 
$H_1$, of $w(M,H_2,~\sim)$.

\begin{theorem}
\label{mwequalsw}
Let $H_1 \subseteq H_2$ be families of subtrees of a given tree.
Let also $\sim$ be a symmetric relation on $H_2$ containing the disjointness
relation (i.e, if $x, y$ are disjoint members of $H_2$ then $x \sim y$).
Then 
$$mw(H_1,H_2,~\sim) = w(H_1,H_2,~\sim)$$ 
\end{theorem}
\vspace{.5cm}

The equality $mw(H) = w(H)$ for a family $H$ of trees is obtained upon
taking $H_1 = H_2 = H$ and taking $\sim$ to be the total relation (i.e., 
$x \sim y$ for all $x, y \in H$). The equality $imw = iw$ is obtained
upon taking $H_1 = H_2 = H$ and $\sim$ the disjointness relation.

The theorem will clearly follow from the following lemma:

\begin{lemma}
Let $H_1, H_2$ and $\sim$ be as above, and let $c,d$ be two intersecting
edges of $H_1$. Then one of $w(H_1 \setminus \{c\},H_2, \sim)$
or $w(H_1 \setminus \{d\},H_2 ,\sim)$ is equal to 
$w(H_1, H_2 ,\sim)$.
\end{lemma}

(Once the lemma is proved, the theorem follows by deleting edges from $H_1$
one by one,
until a matching is obtained).

{\em Proof of the Lemma}
Write $w=w(H_1, H_2 ,\sim)$.
Suppose, for contradiction, that 
$w(H_1 \setminus \{c\},H_2, \sim) = 
w(H_1 \setminus \{d\},H_2, \sim) = w - 1$.
Let $C$ be 
a $\sim$-related covering of size $w-1$
 of $H_1 \setminus \{d\}$
and $D$ a $\sim$-related covering  of size $w-1$ of $H_1 \setminus \{c\}$.
Clearly,
$D$ does not meet $c$ and $C$ does not meet $d$.
Let $S_1$ be the connected part of
$T \setminus (V(c) \cap \bigcup \{V(t): t \in C\})$ containing $d$ and let $S_2 = T \setminus V(S_1)$.
(Here $T \setminus A$, where $A$ is a subset of $V(T)$, denotes the 
forest obtained by deleting the vertices of $A$ from $T$.)

Let 
$C_i = \{g \in C : g \subseteq S_i \}~(i = 1,2)$ 
and
$D_i = \{g \in D : g \subseteq S_1 \}~(i = 1,2)$.
Let also $K = C_1 \cup D_2 \cup \{c\}$ and $L = C_2 \cup D_1$.

Obviously, $K$ is an $H_1$-covering, since the only members of $H_1$
not covered by $C_1$ and $D_2$ are those which meet both $S_1$ and $S_2$,
and these are covered by $c$. 
Since $D$ does
not meet $c$, the tree $c$ is $\sim$-related to all members of $D_2$.
Also, by the definition of $C_1$, no tree from $C_1$ meets $c$, and hence
$c \sim x$ for every $x \in C_1$. Thus $K$ is $\sim$-related.
A similar argument shows that also $L$ is a $\sim$-related $H_1$-covering.
But $|K|+|L| = 2w-1$, and hence either $|K| <w$ or $|L| < w$, 
which contradicts the definition of $w$.

\section{A K\"onig-like theorem}

A hypergraph $H$ is called a {\em point-tree} hypergraph if
it is obtained from a bipartite graph by replacing, in each edge,
the point in one side of the graph (the same side for all edges) 
by a tree. More formally, $H$ is point-tree if,
for some set $X$ and a tree $T$ whose vertex set is disjoint from $X$,
 each edge $e \in H$ is 
of the form $\{x\} \cup V(t)$, where $x =x(e) \in X$ and $t =t(e)$ is
a subtree of $T$. For such a hypergraph we denote by $\sigma(H)$ 
the number $w(H,F)$, where $F = F(H) = \{\{x(e)\}: e \in H\} \cup \{V(t(e)): e \in H\}$.
(That is, the edges of $F$ are to be caught by vertices in $X$, or trees $t(e)$).
Recall that the {\em matching number} $\nu(H)$ of a hypergraph $H$
is the maximal size of a matching in $H$. 
We shall prove  the following generalisation of K\"onig's
theorem:

\begin{theorem}
\label{main}
$\sigma(H) \le \nu(H)$ for every point-tree hypergraph $H$.
\end{theorem}

K\"onig's theorem is the special case in which 
all trees $t(e)$ are singletons.

In \cite{tardos} (see also \cite{kaiser}) there were considered 
objects which are similar
 to point-tree hypergraphs: bipartite graphs in which 
each point (on both sides of the graph) in each edge
is replaced by an interval. In \cite{tardos} it was proved that 
$\tau \le 2\nu$ for such hypergraphs. It is interesting to note that the
proof there was topological (as was the simpler proof in \cite{kaiser}),
and at base so is the proof here, since it relies on Theorem \ref{ah}, whose
proof is topological.

In fact, we shall need a stronger version of Theorem \ref{ah},
its ``deficiency version''.

\begin{definition}
The {\em deficiency} $def({\mathcal A})$ of 
a family 
$\mathcal A$ of hypergraphs is the minimal natural
number $d$ such that $mw(\bigcup{\mathcal B}) \geq |{\mathcal B}|-d$ for every 
subfamily $\mathcal B$ of
${\mathcal A}$.
\end{definition}

The following strengthening of Theorem \ref{ah} 
follows from it by a standard argument (see \cite{ryser}).

\begin{theorem}
\label{deficiency}
Every family $\mathcal A$ of hypergraphs
has a subfamily $\mathcal D$ of size at most $def({\mathcal A})$, 
such that ${\mathcal A} \setminus {\mathcal D}$ has a 
choice function of disjoint edges.
\end{theorem}

{\em Proof of  Theorem \ref{main}}~
For each $x \in X$ let $K(x)$ be the hypergraph 
$\{V(t(e)): x(e) =x\}$.
Applying Theorem \ref{deficiency} to the family $\{K(x): x \in X\}$
yields the existence of a subset $Y$ of $X$
such that, writing $K =
\bigcup\{K(y): y \in Y\}$,
we have $mw(K)= |Y|-(|X| - \nu(H))$. By Lemma \ref{equality} $w(K) = mw(K)$.
Let $C$ be a $K$-covering set of edges in $K$
of size
$|Y|-(|X| - \nu(H))$. 
Then  $C \cup X \setminus Y$
is  an $H$-covering set consisting of elements of $F(H)$,  of size
$\nu(H)$. This proves the theorem. \enp
\\
\
\\

\section{The special case of families of intervals}

A special case of trees is that of intervals. 
Thus,  all the results 
above apply to the special case  of families of intervals. But it turns
out that in this case the equalities $mw = w$ and $imw = iw$
can be strengthened. For
example, the
matchings which are ``hard to catch'' can be described
explicitly, and can be found efficiently. To show this, we shall need the 
following notation. 

We say that a set of vertices in a graph $G$ is {\em $k$-independent}
if the distance between any two of its members is
larger than $k$. (Thus  $1$-independence is just the usual independence.)
By $\gamma_k(G)$ we denote the maximal size of a 
$k$-independent set in $G$.
A set of edges in a hypergraph $H$ is {\em $k$-remote}
if it is $k$-independent in the line graph $L(H)$.
We write $\zeta_k(H)$
for $\gamma_k(L(H))$. Clearly, for all hypergraphs we have

\begin{equation}
\label{ineq}
\zeta_2 \leq mw \leq w
\end{equation}

Here we shall prove:

\begin{lemma}
\label{equality}
For a family $F$ of 
intervals on the real line  equality holds throughout  (\ref{ineq}),
namely
$\zeta_2(F) = mw(F)= w(F)$.
\end{lemma}

Given a family of intervals $F$, we 
denote by $\ell(F)$ the leftmost interval in $F$, namely the one 
with minimal right endpoint.
For a real number $x$ we denote by $F(>x)$ the set of those intervals
in $F$ which lie entirely within $(x, \infty)$.

{\em Proof of Lemma \ref{equality}}~
By (\ref{ineq}) the lemma will be proved if we can exhibit
a $2$-remote set $R$ and a covering set $C$ of the same size.
We construct such sets by an inductive process.

Let $r_1 = \ell(F)$, and 
$c_1 = [x_1,y_1]$ an element of $F$ meeting $r_1$ whose right
endpoint is largest. 
Let now $r_2 = \ell(F(>y_1))$,
and 
$c_2 = [x_2,y_2]$   the interval in $F$ meeting $r_2$
whose right endpoint is maximal. (Note that $c_2$ may intersect $c_1$.) 
In general, having defined $c_i=[x_i,y_i]~(i>0)$,
we let
$r_{i+1} = \ell(F(>y_i))$ and 
choose $c_{i+1}$ as the interval in $F$ meeting $r_{i+1}$ and extending
furthest to the right.
The process terminates when, at some stage $t$,
there is no interval lying entirely to the right of $c_t$. Then, clearly, 
$R=\{r_1, \ldots, r_t\}$ is remote, while $C=\{c_1,\ldots, c_t\}$
is covering.
\enp

A similar equality can be proved for $\zeta_k$ for all $k$.
For any graph $G$ we denote by $\rho_k(G)$ the minimal number
of vertex sets  of diameter at most $k$, whose 
union is $V(G)$. Obviously, $\rho_k(G) \ge \gamma_k(G)$ for all $k$.
One corollary of Lemma \ref{equality} is that 
in an interval graph
$\rho_2(G) = \gamma_2(G)$ (in fact, it implies something even  
stronger, namely 
that an interval graph $G$ can be partitioned into 
$\gamma_2(G)$ sets of radius $1$.)
This equality is true for all $k$, that is, 
in an interval graph
$\rho_k(G) = \gamma_k(G)$ 
for all $k$. We shall prove this for a more general 
class of graphs, that of {\em incomparability graphs}.
A graph is called an {\em incomparability graph}
if there exists a partial order on its vertex set, such that 
two vertices in the graph are connected if and only if they are 
incomparable. As is well known, an interval graph is 
an incomparability graph, the order on the intervals being that of
``lying completely to the right of''.

\begin{proposition}
\label{rhogamma}
In an incomparability graph 
$\rho_k(G) = \gamma_k(G)$ for all $k$.
\end{proposition}

Given a graph $G$ and a natural number $k$, we denote by $G^{*k}$ the graph 
with the same vertex set as $G$, in which two vertices are connected if and
only if their distance in $G$ does not exceed $k$. Obviously, 
$\gamma_k(G) = \gamma(G^{*k})$, and $\rho_k(G) = \rho(G^{*k})$. Thus, 
if $G^{*k}$ is perfect, then $G$ satisfies Proposition \ref{rhogamma}.
Since incomparability graphs are known to be perfect, 
Proposition \ref{rhogamma} will follow from:

\begin{lemma}
If $G$ is an incomparability graph 
then so is \gs.
\end{lemma}

\pf
Let ``\ger'' be a partial order on $V(G)$ such that 
$x$ and $y$ are connected in $G$ if and only if 
they are \ger-incomparable.
Let ``\gep'' be the relation defined by 
$x \gep y$ if and only if $x \ger y$ and
$(x,y) \not \in E(\gs)$. Clearly, it suffices to show that 
``\gep'' is a partial order. 

Assume that it is not. Then there exist $x, y, z \in V(G)$
such that $x \gep y, ~y \gep z$, while $x \not \gep z$.
There exist then vertices $x = u_0, u_1, u_2, \ldots ,u_t=z~(t \le k)$
such that $u_i$ is $\ger$-incomparable to $u_{i+1}$
for all $0 \le i < t$. Since $x \gep y$, all $u_i$ must be
$\ger$-comparable to $y$ (or else 
the sequence $u_0, \ldots, u_i, y$ would show that 
$(x,y) \in E(\gs)$). Since $u_0 \ger y$ and $y \ger u_t$,
it follows that there exists  $i < t$ 
such that $u_i \ger y$ and $y \ger u_{i+1}$. But then, by 
the transitivity of the relation $\ger$,
we have $u_i \ger u_{i+1}$, a contradiction.
\enp

Our next observation on the interval case is that 
the equality $imw = iw$ can be given in this case a 
constructive (and algorithmically efficient) proof.

\begin{theorem}
\label{indmatching}
$imw(F) = iw(F)$ for every  hypergraph $F$ of intervals.
Moreover, there is a polynomial time algorithm producing a 
matching $M$ such that $iw(M,F) = iw(F)$.
\end{theorem}

{\em Proof~}
Write $k = iw(F)$. We shall prove the theorem by constructing
a matching
$M$ in $F$ satisfying:
\\
\
\\
(P)~~ Every matching of size $k-1$ in $F$
misses an interval from $M$. 
\\
\

Note, first, that we can assume that no two endpoints of intervals in $F$
coincide: if there is such coincidence, we can enlarge some intervals,
without changing the intersection pattern of the family. Henceforth we shall
make this assumption.

For a matching $E$ in $F$ we write 
$E = (e_1, \ldots, e_t)$ if $e_i$ 
are the edges of $E$ ordered from left to right.
Let $E$ be such a matching, and write
$e_j = [a_j, b_j], ~j \le t$. 
For each $j \le t$  we denote by $d_j(E)$ the interval
$\ell(F(>b_{j-1}))$ 
(where $b_0$ is defined as $-\infty$, meaning
that $d_1(E)$ is always $\ell(F)$).
If $F(>b_t)$ is non-empty,
we write also $d_{t+1}(E)=\ell(F(>b_t))$.
Note that $d_j(E)$ meets $e_j$ if and only if 
there is no interval from $F$ lying between 
$e_{j-1}$ and $e_j$; or, in the case $j=1$, if and only if there
is no interval in $F$ lying to the left of $e_1$.
By $m(E)$ we denote the 
maximal index $m \le t$ such that 
there is no interval from $F$ lying between 
$e_{j-1}$ and $e_j$
for all $j \le m$ (if there is an interval 
from $F$ to the left of $e_1$, we write $m(E) = 1$).
If $m(E) = t$, we say that $E$ is {\em dense}.
Note that if $E$ misses some edge from $F$ then it misses
$d_{m(E)}(E)$.

Let $L_i$ be the set of all right endpoints 
of intervals $e_i$ in some 
dense matching $E = (e_1, \ldots, e_t)$.
An interval $[a,b]$ is called $L_i$-{\em free} if
$[a,b] \cap L_i \subseteq \{b\}$.
Let $M$
be the set of all intervals in $F$ which are 
equal to $d_i(E)$ for some dense matching $E$, and are
$L_j$-free for all
$j \le i$.
\\
\
The proof will be complete if we show that $M$ is a matching, and that 
it satisfies property (P).
\\

Let us first show that $M$ is a matching. Suppose, for contradiction,
that $d = d_p(E) = [a,b]$ and $d' = d_q(E') = [a',b']$  
are two intersecting members of $M$,
where $E=(e_1,\ldots ,e_{p-1})$
and $E'=(e'_1,\ldots ,e'_{q-1})$
are two dense matchings. By assumption
$b \neq b'$, so we may assume (say) that
$b < b'$. Note that $b \in L_p$, since 
the matching $(e_1, \ldots e_{p-1},d)$ is dense.
Since
$d' \in M$, this implies that 
$q<p$. Let $y$ be the right endpoint of $e'_{q-1}$.
Since $d \in M$ and $y \in L_{p-1}$, it follows that 
$y < a$. But then the interval $d=[a,b]$ witnesses
the fact that $d' \neq \ell(F(>y))$, contradicting the 
definition of $d'$.

Next we show that $M$ satisfies property (P).
Suppose that there exists a matching $Z= (z_j: j \le t)$ 
in $F$ meeting all
intervals in $M$, where $t < k$. 
Among all such matchings $Z$, choose one 
in which the right endpoint of $z_{m(Z)}$
is maximal. Write $m=m(Z)$ for this $Z$.
Since $iw(F) = k$, the matching
$Z$ misses an interval from $F$. This means that 
$Z$ misses $d = d_m(Z)$. If $d \in M$, we are done. If not, then
$d$ contains a point from $L_j$ for some $j \le m$.
That is, it is the right
endpoint of an interval $w_j$ in some dense matching 
$W=(w_1, \ldots ,w_m)$. But then the matching
$(w_1, \ldots, w_j, z_{m+1}, \ldots ,z_t)$, like $Z$,
does not miss any edge from $M$. This contradicts the maximality 
property of $Z$.

To see that $M$ above can be constructed in polynomial time,
it suffices to note that the sets $D_j$ of those intervals which appear as
the $j$-th element in some dense matching can be constructed inductively in
plynomial time, and that $M$ is defined by these sets.
\enp

\end{document}